\documentclass[11pt,leqno]{article}
\usepackage{amsmath, amsthm, amsfonts}
\begin{document}
\newtheorem{theo}{Theorem}[section]
\newtheorem{defi}{Definition}[section]
\newtheorem{cor}{Corollary}[section]
\newtheorem{lem}{Lemma}[section]
\newtheorem{pro}{Proposition}[section]
\newtheorem{que}{Question}[section]
\newtheorem{rem}{Remark}[section]
\newtheorem{exa}{Example}[section]
\newtheorem{conj}{Conjecture}[section]
\def\pr{{\bf Proof: }}
\def\Sff{{\cal F}(d)}
\def\M{{\cal M}  (d)}
\def\Mz{{\cal M } _0 (d)}
\def\proj{{\Bbb P}^2} 
\def\pl{{{\Bbb P}^1}}
\def\integ{\frac{F^p}{G^q}}
\def\O{{\cal O}} 
\def\C{{\Bbb C}}
\def\S{{\cal S}}
\def\Z{{\Bbb Z}}
\def\P{{\cal P}}
\def\Ci{{\cal C}}
\def\R{{\Bbb R}}
\def\RR{{\cal R}}
\def\F{{\cal F}}
\def\Sf{{\cal F}(M,L)}
\def\Hf{H^0(M,\Omega^1\otimes L)}
\newcommand{\cfv}[1]{\{\delta_t\}_{t\in #1}}
\newcommand{\form}[1]{\Omega^1(#1D)}
\def\Pf{\C[t]}
\begin{center}
{\LARGE\bf Abelian integrals in holomorphic foliations\\}
\vspace{.25in} {\large {\sc Hossein Movasati}\footnote{
Keywords: Holomorphic foliations - holonomy - Picard-Lefschetz theory
\\
Math. classification: 57R30 - 14D99 - 32G34}}
\\
{\it Dedicated to Mothers}
\end{center}
\begin{abstract}
The aim of this paper is to introduce the theory  of Abelian
integrals for holomorphic foliations in a complex manifold of
dimension two. We will show the importance of Picard-Lefschetz
theory  and the classification of relatively exact 1-forms in
this theory. As an application we identify some
irreducible components of the space of holomorphic foliations of
a fixed degree and with a center singularity in the projective space of
dimension two. Also we  calculate higher Melnikov functions under
some generic conditions.
\end{abstract}
\setcounter{section}{-1}
\section{Introduction}
Let us be given a differential equation
\begin{equation}
\label{14jan01}
\frac{dy}{dx}=\frac{P(x,y)}{Q(x,y)}
\end{equation}
in the real plane ${\Bbb R}^2$, where $P$ and $Q$ are two
polynomials in $x$ and $y$. Let $H(P,Q)$ denote the number of
limit cycles of the above differential equation and
$$
H_n=max\{H(P,Q)\mid deg(P),deg(Q)\leq n\}
$$
For many years it was believed that the finiteness of $H(P,Q)$ have been 
proved by Dulac in \cite{dul}. But it was found that \cite{dul} contains
a gap and then the complete proof was given in \cite{ily1} and \cite{eca} 
independently. Nowadays the question of whether $H_n$ is finite or not
is considered as the Hilbert sixteen problem in differential equations.
 An equation ~(\ref{14jan01})
 has a first integral or is called
integrable if there are two polynomials $F$ and $G$ in ${\Bbb R}^2$ such that $\frac{F}{G}$
is constant on its solutions. In this case the equation (after reducing) has the form
$$
\frac{dy}{dx}=-\frac{GF_x-FG_x}{GF_y-FG_y}
$$
When $G\equiv 1$ then the equation is called Hamiltonian. 
Let $\S'_n$ (resp. $\S_n$) denote the space of equations
~(\ref{14jan01}) defined by real (resp. complex) polynomials $P$
and $Q$ $deg(P),deg(Q)\leq n$. The space $\S'_n$ is
parameterized  by the coefficients of $P$ and $Q$.

One of the first attempts to solve the mentioned problem was
made  by two Russian mathematicians I. G. Petrovski\u\i\ and E.
M. Landis (\cite{lape1},\cite{lape2}). They complexified the
equation ~(\ref{14jan01}) and considered the equation in $\C^2$
 and tried to find another set of cycles $\Ci(P,Q)$ in the solutions 
of ~(\ref{14jan01}) whose  cardinality is not less than
 the number of limit cycles in the related
real equation. Then they tried to calculate a uniform upper bound
for the cardinality of ${\cal C}(P,Q)$ (and hence $H_n$) by
perturbation of integrable equations. 
Let $\{\delta_t\}_{t\in(\R,0)}$ be a continuous family of closed 
solutions of an
integrable equation. They observed that the bifurcation of limit
cycles from this family of closed cycles is related with the
zeros of a certain Abelian integrals.

Although their work had errors, the idea of using Abelian
integrals to count the number of limit cycles of a perturbed
Hamiltonian equation has been one of the useful tools in getting
approach to the Hilbert Sixteen Problem. For this see the book
\cite{rou} and its references.
The aim of this text is to introduce an algebro-geometric
approach to these Abelian integrals. Instead of $\C^2$ we consider
an arbitrary two dimensional compact complex manifolds and instead
of Hamiltonian fibrations we consider the fibrations of
meromorphic functions on $M$. In section ~\ref{foliation} and
~\ref{piclef} holomorphic foliations and the Picard-Lefschetz
theory of meromorphic functions on $M$ are introduced. In section
~\ref{akharinseminar} we see how Abelian integrals appear in the
formula  of the first Melnikov function. In section
~\ref{30jul01} we will encounter with relatively exact 1-forms.
In section \ref{hargez} we will classify some irreducible components
of holomorphic foliations in $\proj$ with a center singularity
and finally in section ~\ref{melnikov} we will calculate higher
order Melnikov functions under some generic conditions.

This work started from my doctorate thesis. I want to use this opportunity to express my thanks to my advisor Alcides Lins Neto, who introduced me with the article \cite{ily}. Also, thanks go to Cesar Camacho and  Paulo Sad
my teachers in IMPA-Brazil.
I thank also S. Shahshahani who made my visiting from IPM-Iran possible
when I was working on some parts of this work. 
I thank also MPIM-Germany which the final version of this paper was obtained there.
\section{Holomorphic Foliations}
\label{foliation} Let $M$ be a compact complex manifold of
dimension two. A holomorphic foliation in $M$ with isolated
singularities is given by a collection of holomorphic 1-forms
$\omega_\alpha$ defined on $U_\alpha, \alpha\in I$, where
$\{U_\alpha\}_{\alpha\in I}$ is an open covering of $M$, and such
that
\begin{equation}
\label{shabjan}
\omega_\alpha=g_{\alpha\beta}\omega_\beta,\  \alpha,\beta\in I
\end{equation}
where $g_{\alpha\beta}$ is a  holomorphic without zero function
in $U_\alpha\cap U_\beta$. Furthermore we assume that the set of
points in which $\omega_\alpha$ is zero has codimension greater
than one (discrete set). In other words $\omega_\alpha$ has not
zero divisor.

Therefore for any foliation $\F$ there is associated a line
bundle $L$ given by the transition functions
$\{g_{\alpha\beta}\}_{\alpha,\beta\in I}\in H^1(M, \O^*)=Pic(M)$.
We also say that the foliation $\F$ is of degree $L$. The data
~(\ref{shabjan}) can be considered as a holomorphic section
$\omega\in \Hf$ without zero divisor, where $\Omega^1$ is the
cotangent bundle of $M$. Since $M$ is compact, $\Hf$ is a finite
dimensional space (this is a simple corollary of Grauert direct
image theorem). By $\F(\omega)$ we mean that the foliation $\F$ is
given by the 1-form $\omega\in \Hf$. It is easy to check that for
any foliation $\F$ there exists a unique line bundle $L$ in $M$
such that $\F$ is given by a holomorphic without zero divisor
element of $\Hf$. From now on we fix a line bundle $L$ and assume
that $\Hf$ has a holomorphic without zero divisor section. In this
case the set of 1-forms $\omega\in \Hf$ without zero divisor is
an open subset of $\Hf$.

If $\omega\in \Hf$ has a zero divisor we use the following trick:
Let $\omega\in \Hf$ be a holomorphic section with the zero
divisor $Z$. Let $L_Z$ be the line bundle associated to $Z$ and
$s\in H^0(M, L_Z)$ be the holomorphic section with the zero
divisor $Z$. Now $\frac{\omega}{s}$ is a holomorphic without zero
divisor  section of $H^0(M, \Omega\wedge L\wedge L_Z^{-1})$, and
so, we can substitute $L$ for $L\wedge L_Z^{-1}$.

Two holomorphic  without zero divisor sections $\omega,\omega'\in
\Hf$ induce the same foliation if and only if
$\omega=c.\omega'$, where $c$ is a constant. Therefore the space
of foliations of degree $L$, namely $\F(M,L)$, form an open subset
of the projective space
$$
Proj(\Hf)
$$
Many times we need the foliation $\F$ to be given by a
meromorphic 1-form in $M$. In these cases we assume that $H^0(M,
L)\not=0$ and choose a non-zero section $s$ of $H^0(M, L)$. Now
the foliation $\F(\omega),\omega\in \Hf$ is given by the
meromorphic 1-form $\frac{\omega}{s}$ in $M$.
\begin{defi}\rm The
meromorphic section $s\in H^0(M, L)$ is called the integrating
factor of $\F(\omega),\omega\in \Hf$ if $\frac{\omega}{s}$ is a
closed meromorphic function in $M$. $\F(\omega)$ is called
integrable if there exists an integrating factor $s\in H^0(M, L)$
such that $\frac{\omega}{s}=df$, where $f$ is a meromorphic
function on $M$. In this case $f$ is constant on the leaves of
$\F$. We also say that $f$ is a first integral of $\F$.
\end{defi}
Let $\F$ be an integrable  holomorphic foliation  given by a
without zero holomorphic 1-form $\omega\in\Hf$. Let also $f$ be a
first integral of $\F$. The fibers of $f$ has the same linear
bundle which we denote it by $L_A$, where $A$ is a generic fiber
of $f$ . Let $A_1,A_2,\ldots,A_s$ be irreducible components of
the fibers of $f$ in which $f$ has multiplicities greater than
one, namely $n_1,n_2,\ldots,n_s$, respectively. It is easy to see
that
$$
 L= L_{2A-\sum_1^s(n_i-1)A_i}
$$
For the projective space $\proj$ we have $Pic(\proj)\cong\Z$. We can take
a line $H$ as the generator of $Pic(\proj)$.  Therefore every line bundle
in $\proj$ is of the form $L_{d.H}$, where $d$ is an integer. It is easy to 
see that for a line bundle coming from a foliation $d$ is positive. In this case $d$ is called the degree of the foliation. Every degree $d$ holomorphic foliation in $\proj$
is given by a polynomial 1-form $Pdx-Qdy$ in the affine chart $\C^2=\proj-H$.
\section{Picard-Lefschetz Theory}
\label{piclef}
\def\Q{{\Bbb Q}}
\label{1feb01} This section is mainly based on the papers
\cite{lam} and \cite{hos2}. For more information the reader is
referred to that papers. Throughout the text when we do not write
the coefficients used in the homology we mean homology with
coefficients in $\Z$.
\\
Let $f$ be a meromorphic function on a complex manifold $M$
of dimension two. The indeterminacy set $\RR$ of $f$ contains the points of $M$ in which $f$ has the form $\frac{0}{0}$. $\RR$ is a discrete set and  the following
holomorphic function is well-defined:
$$
f:M-\RR\rightarrow\pl
$$
We use the following notations
$$
L_K=f^{-1}(K),\  M_K=\overline{L_K}, \ K\subset \pl
$$
For any point $c\in \pl$ by $L_c$ and $M_c$ we mean the set
$L_{\{c\}}$ and $M_{\{c\}}$, respectively. Throughout the text by
a compact $f$-fiber we mean $M_t$ and by a $f$-fiber only we mean
$L_t$. Performing a finite number of blow-ups (see \cite{cam}) at
the points of $\RR$ and using Ehresmann's Fibration Theorem (see
\cite{lam}), we can see that there exists a finite subset
$C=\{c_1,c_2,\ldots,c_r\}$ of $\pl$ such that $f$ fibers $M-\RR$
locally trivially over $B=\pl-C$ i.e., for every point $b\in B$
there is a neighborhood $U$ of $b$ and a $C^\infty$-diffeomorphism
 $\phi :U\times f^{-1}(c) \rightarrow f^{-1}(U)$ such that
 $f\circ\phi =\pi _1 =$ the first projection. We say that $C$ is the set of critical values
 of $f$.
 The regularity of $f$ along a fiber $L_c$ does not imply that it is fiber bundle over a
 neighborhood of $c$. This situation happens when $M_c$ has
 a tangency point with the divisor of the blow-up in which the
 leaves separate from each other. In this case we say that the
 critical point related to the value $c$ is in $\RR$.

Now we are able to use the Picard-Lefschetz theory for
understanding the topology of the fibers of $f$.
 Let $\lambda$ be a path in $B$ with the initial
and end points $b_0$ and $b_1$.
There is an isotopy
$$
H:L_{b_0}\times[0,1]\rightarrow L_{\lambda}
$$
such that for all $x\in L_{b_0},\ t\in [0,1]$
$$
 H(x,0)=x,\ H(x,t)\in L_{\lambda (t)}
$$
For every $t\in [0,1]$ the map $h_t=H(.,t)$ is a homeomorphism between
$L_{b_0}$ and $L_{\lambda (t)}$. The different choices of $H$
 and paths homotopic
 to $\lambda$ in $B$ would give the class of isotopic maps
$$
 \{h_{\lambda} : L_{b_0} \rightarrow L_{b_1}\}
$$
 where $h_\lambda(.)=H(.,1)$. The class $\{h_{\lambda} : L_{b_0} \rightarrow L_{b_1}\} $
 defines the map
$$
h _{\lambda} :H_1(L_{b_0}) \rightarrow H_1(L_{b_1})
$$
For any regular value $b$ of $f$, we can define
$$
h:\pi _1(B,b)\times H_1(L_b)\rightarrow H_1(L_{b})
$$
$$
h(\lambda ,.)=h_\lambda(.)
$$
The action of $\pi_1(B,b)$ through $h$ on $H_1(L_b)$ is called the action of monodromy on the first
 homology group of
$L_b$. We also say that $H_1(L_c)$ is a $\pi_1(B,b)$-module.

Let $q$ be an indeterminacy point of $f$. Let $\tau: \tilde{M}\rightarrow
M$ be the map of successive blow-ups in $q$ such that $f\circ\tau$ is 
holomorphic in $\tau^{-1}(M,q)$ and there is an irreducible component
$D_1\cong\pl$ of the divisor $\tau^{-1}(q)$ such that $f\circ\tau\mid_{D_1}$ is 
an isomorphism (see \cite{cam} for details of such blow-ups). 
Every $\tau^{-1}(M_t)$ passes through $D_1$ transversally. Therefore we can
choose a continuous family of small cycles in $\tau^{-1}(M_t)$ around $D_1$.
The image of this family by $\tau$ is simply called 
a continuous family of cycles around $q$. This family is not defined for
$t\in \pl$ such that $\tau^{-1}(M_t)$ contains an irreducible 
component of $\tau^{-1}(q)$. Particularly it is defined for all $t\in\pl-C$. 
In fact a $\delta_t$ is obtained by
monodromy along an arbitrary path from a fixed $\delta_{b}$. The action of monodromy on $\delta_{b}$ is identity. 

\begin{defi}\rm
Let $K$ be a subset of $B$ and  $b$ be a point in $ K\backslash C$.
Any relative 2-cycle of $L_K$ modulo $L_b$ is called  a 2-thimble above
$(K,b)$ and its boundary in $L_b$ is called a vanishing 1-cycle above $K$.
\end{defi}
Suppose that $f$ has a non-degenerate critical point at $p_i\in
M-\RR$ and $f(p_i)=c_i$. Let $\lambda _i$ be  a path in $B\cup
\{c_i\}$ connecting a regular value $b$ to $c_i$. This path
passes through $c_i$ only in its end point. The leaf $L_b$
contains a cycle which vanishes above $\lambda_i$ and called the
Lefschetz vanishing cycle. Roughly speaking, when the value $x$
varies from $b$ to $c_i$, the Lefschetz vanishing cycle in $L_x$
moves in the leaves and arrives to the leaf $L_{c_i}$ as the
point $p_i$. The locus of this movement is exactly the related
2-thimble.

Fix a point $\infty\in\pl$, which may be a critical value.
Let $C'=\{c_1,c_2,\ldots,c_r\}$ be the subset of $\C=\pl-\{\infty\}$ 
containing critical
values of $f$.

Consider a system of $r$  paths $\lambda _1,\ldots,\lambda _r$ starting
from $b$ and ending at $c_1, c_2,\ldots , c_r$, respectively, and such that:
\begin{enumerate}
\item
 each path $\lambda _i$ has no self intersection points  ;
\item
two distinct path $\lambda _i$ and $\lambda _j$ meet only at their
 common origin $\lambda _i(0)=\lambda _j(0)=b$.
\end{enumerate}
This system of paths is called  a distinguished system of paths.
Let $K$ be the union of these paths and small disks around
$c_i$'s.
 The set of vanishing cycles above $K$ in $L_b$ is called a distinguished set of vanishing
 cycles related to the critical points $c_1,c_2,\ldots,c_r$.

\begin{theo}(Theorem 2.2.1 \cite{hos2})
\label{31dec01}
Suppose that $H_{1}(M-M_\infty,\Q)=0$.
 Then a distinguished set of vanishing 1-cycles
 related to the critical points in the set
  $C\backslash \{\infty\}=\{c_1,c_2,\ldots,c_r\}$ generates $H_1(L_b,\Q)$.
\end{theo}
Again note that in the above theorem $\infty$ can be a critical
value of $f$. If $M_\infty$ is smooth and its homology class in
$H_2(M,\Q)$ is not zero (for example if $M_\infty$ is a smooth
hyperplane section) then the condition $H_{1}(M,\Q)=0$ implies
that $H_{1}(M-M_\infty,\Q)=0$. The reason is as follows: By Leray
(or Thom-Gysin) isomorphism (see \cite{che} p. 537) we have
$H_2(M,M-M_\infty)\simeq H_0(M_\infty)\simeq \Z$. Now we write
the long exact sequence of the pair $(M,M-M_\infty)$
$$
\cdots\rightarrow H_2(M-M_\infty)\stackrel{i}\rightarrow
H_2(M)\rightarrow \Z\stackrel{j}{\rightarrow}
H_1(M-M_\infty)\rightarrow 0
$$
Therefore $H_1(M-M_\infty)=\Z/n\Z$ or $\Z$ for some natural number
$n$. In the first case we have $H_1(M-M_\infty,\Q)=0$. In the
second case $j$ is one to one and so $i$ is surjective. But this
means that the intersection of $M_\infty$ with any 2-cycle in
$H_2(M)$ is zero. By Poincar\'e duality the class of $M_\infty$ in
$H_2(M,\Q)$ must be zero which is a contradiction.

\begin{defi}
\label{11jan02}
The cycle $\delta$ in a regular fiber  $L_b$ is called simple if
the action of $\pi_1(B,b)$ on $\delta$ generates $H_1(L_b,\Q)$.
\end{defi}
Note that in the above definition we have considered the homology group with 
rational coefficients. Of course, not all cycles are simple. 
For instance if the meromorphic function in a local coordinate $(x,y)$ 
around $q\in \RR$ has the form $\frac{x}{y}$, then the cycle around $q$ 
in each leaf has this property that it is fixed under the action of 
monodromy, therefore it cannot be simple. In the next paragraph we 
are going to introduce some fibrations with simple cycles.

\def\pn{{\Bbb P}^n}
\def\pnh{\check{{\Bbb P}}^n}
{\bf Lefschetz pencil: }The hyperplanes of $\pn$ are  points of the dual projective space $\pnh$.
We use the following notation:
$$
H_y \subset \pn\ , \ y\in\pnh
$$
Let $M$ be a closed irreducible smooth subvariety of $\pn$. Its dual variety $\check{M}$ consists of all points $y$ in $\pnh$ such that $H_y$ is tangent to $M$ at some point. $\check{M}$ may have singularities.
A Lefschetz pencil in $\pn$ consists of  all hyperplanes which contain a fixed
$(n-2)$-dimensional projective space  $A$, which is called the axis
of the pencil. We denote a pencil by
$$
\{H_t\}_{t\in G}
$$
or $G$ itself, where $G$ is a projective line in $\pnh$.
The pencil $\{H_t\}_{t\in G}$ is in general position with respect to $M$
 if $G$ is in general position with respect to $\check{M}$.
Sometimes we parameterize the line $G\cong \pl$ and consider the meromorphic function $f$ on $M$ induced by the pencil.
 The meromorphic function associated to a pencil in general position is called
 the generic Lefschetz meromorphic function of the pencil 
and has the following properties:

\begin{enumerate}
\item
In a local coordinates system $(x_1,x_2,\ldots,x_{n})$ around an 
indeterminacy point $q$ of $f$, $f$ has the
form $\frac{x_1}{x_2}$;
 \item
$f$ has $r=deg(\check{M})$ non-degenerate
critical points $p_1,\ldots, p_r$ in $M-A$ such that $f(p_i)=c_i$'s
are distinct values in $\pl$;
\item
For every critical value $c_i$ of $f$, the hyperplane
 $H_{c_i}$ has a unique tangency of order two with $M$ which lies out of $A$.
 The other hyperplanes  are transverse to $M$
 (see \cite{lam}).
\end{enumerate}

\begin{theo}
\label{30jan01} Suppose that $H_1(M-M_\infty,\Q)=0$ and
$H_1(M_\infty)\not = 0$, where $M_\infty$ is a regular fiber of a
generic Lefschetz meromorphic function $f$. Then every Lefschetz
vanishing cycle in a regular fiber of $f$ is simple.
\end{theo}
Note that the regular compact fibers of $f$  have the same
topology and so $H_1(M_\infty)\not =0$ means that the genus of
any regular compact fiber is  greater than 0 and so the fibration
is not rational.
\\
Almost all the arguments to prove the above theorem exist in
\cite{lam}. That article has worked with the topology of $M_t$'s,
and not $L_t$'s. Therefore we have proved above theorem in
\cite{hos2} Theorem 2.3.2.

Let $\F$  be a foliation in $\proj$ with a first integral of the
type $\integ$, where $F$ and $G$ are  two relatively prime
irreducible polynomials in an affine chart $\C ^2$ of $\proj$,
$\frac{deg(F)}{deg(G)}=\frac{q}{p}$ and $g.c.d.(p,q)=1$. Assume
that $\F$ has the following generic properties:
\begin{enumerate}
\item
 $\{F=0\}$ and $\{G=0\}$ are smooth and intersect each other transversally;
\item
The critical points of $\integ$ in $\proj\backslash(\{F=0\}\cup\{G=0\})$ are
non-degenerate with distinct images.
\end{enumerate}
 Let $c_1,c_2,\ldots,c_r$ denote these critical values.  In \cite{hos2} it
 is proved that these conditions are generic, i.e.  there exists a dense open
  subsets in the space of coefficients of $F$ and $G$ such that all polynomials $F$ and
  $G$ whose coefficients are chosen from this set, satisfy the
conditions 1 and 2.
The meromorphic function $\integ$ is a fiber bundle over
$\pl-\{c_1,c_2,\ldots,c_r,0,\infty\}$. It has multiplicity $p$ along $F=0$ and $q$ along $G=0$.
\begin{theo}
(\cite{hos2}) \label{argentina} If $deg(F)+deg(G)>4$ then every
Lefschetz vanishing cycle in a regular fiber of $\integ$ is
simple.
\end{theo}
It is easy to see that $\integ$ is a rational fibration (a
fibration whose generic fiber
 has genus zero) if and only if $deg(F)+deg(G)\leq 2$. 
In fact only in the case 
$deg(F)+deg(G)\leq 2$ the
 fibration $\F$ has reducible critical fibers. In the case 
$def(F)=deg(G)=1$ Theorem ~\ref{argentina}
 is not true. In this case it is easy to
 find counterexamples for the forthcoming Theorem ~\ref{main} by performing deformations inside the logarithmic foliations.

 If
$p=q=1$ then $\integ$ is a Lefschetz meromorphic function and the
above theorem is a particular case of Theorem ~\ref{30jan01}. In
the case where $G$ is a linear polynomial we have $p=1$ and
$q=deg(F)$. Considering $G=0$ as the line at infinity in $\proj$,
we have the Hamiltonian fibrations in $\C^2$. In this case the
above theorem is proved by Ilyashenko \cite{ily} using a theorem
of Zizcenko.
\section{Deformation of Holomorphic Foliations and Abelian
Integrals} \label{akharinseminar} Let $\F$ be an integrable
holomorphic foliation in $\Sf$  and
$\F_\epsilon(\omega_\epsilon)\in\Sf, \epsilon\in (\C,0)$  a
holomorphic deformation of $\F=\F_0$. The set $\Hf$ is a vector
space and so we can write
$$
\omega_\epsilon=\omega_0+\epsilon\omega_1+\epsilon^2\omega_2+\cdots,\ \omega_i\in \Hf
$$
The 1-form $\omega_1$ is called the tangent vector of the deformation.
\\
Let $\delta$ be a cycle in a leaf of $\F$ and $\Sigma\simeq(\C,0)$
a holomorphic section to $\F$ in a point $p\in \delta$. Let also
$s$ be an integrating factor of $\omega_0$ whose zero divisor
does not intersect $\delta$ and
$$
\frac{\omega_0}{s}=df
$$
where $f$ is a meromorphic function on $M$.

Throughout the text we  assume that the transverse section
$\Sigma$ is parameterized by $t=f\mid_\Sigma$.
 Assume that the holonomy of $\F$ along $\delta$ is identity.
Note that if $f$ has multiplicity one along the leaf containing $\delta$ then the holonomy is always identity.
We can consider the collection of $\F_\epsilon$'s  as a codimension 
two foliation $\bar {\F} =\{\F _\epsilon\}_{\epsilon\in (\C,0)}$ in
$ M \times (\C ,0)$ and $\Sigma\times (\C ,0)$ as a transverse section to $\bar{\F}$. So we have the holonomy map defined by
$$
H:\Sigma\times (\C,0)\rightarrow \Sigma\times(\C,0)
$$
$$
H(t,\epsilon)=(h_\epsilon(t),\epsilon)
$$
$h_\epsilon(t)$ is a holomorphic function in $\epsilon$ and $t$ and  is
called the holonomy
of $\F_\epsilon$ along the path $\delta$ (Note that by hypothesis $h_0(t)\equiv t$).
 We write
$$
h_\epsilon(t)-t=M_1(t)\epsilon+M_2(t)\epsilon^2+\cdots +M_i(t)\epsilon^i+\cdots,\  i!.M_i(t)=\frac{\partial^ih_\epsilon}{\partial \epsilon^i}\mid_{\epsilon=0}
$$
$M_i$ is called the $i$-th Melnikov function of the deformation along the path $\delta$. Let $M_1\equiv M_2\equiv\cdots\equiv M_{k-1}\equiv 0$ and $M_k\not\equiv 0$. It is a well known fact that the multiplicity of $M_k$ at $t=0$ is the number of limit cycles (more precisely the number of fixed points
of the holonomy $h_\epsilon $) 
which appears around $\delta$ after the deformation. This fact shows the importance of these functions in the local study of Hilbert 16-th problem. The following proposition gives us a nice formula for the first Melnikov function.
\begin{pro}
\label{1febb}
The first Melnikov function is given by
$$
M_1(t)=-\int_{\delta_t}\frac{\omega_1}{s}
$$
where $\omega_1$ is the tangent vector of the deformation and $\delta_t$ is the lifting up of $\delta$ in the leaf through $t\in \Sigma$ .
\end{pro}
\pr The proof is completely formal in the literature of
differential equations (see \cite{rou} and \cite{fra}). The
deformed foliation is given by the meromorphic 1-form
\begin{equation}
\label{sumiu}
df+\epsilon\frac{\omega_1}{s}+O(\epsilon^2)
\end{equation}
Let  $\delta_{t,h_\epsilon(t)}$ be a path in the leaf of $\F_\epsilon$
through $t$ which connects $t$ to $h_\epsilon(t)$ along the path $\delta$.
Since $\Sigma$ is parameterized by $t=f\mid_\Sigma$, by integrating the 1-form ~(\ref{sumiu}) over the path $\delta_{t,h_\epsilon(t)}$ we have
$$
h_\epsilon(t)-t+\epsilon(
\int_{\delta_t}\frac{\omega_1}{s}+O(\epsilon))+O(\epsilon^2)=0
$$
The coefficient of $\epsilon$ in the above equality gives us the desired equality.
\qed

We want to have an explicit formulas for higher Melnikov functions. For this purpose we must classify a certain class of relatively exact 1-forms modulo an integrable foliation.
\section{Relatively Exact 1-forms} \label{30jul01}
\def\L{{\cal L}}
First let us give the definition of a relatively exact 1-forms modulo a foliation.
\begin{defi}\rm
Let $\F$ be a foliation in $M$. A meromorphic 1-form $\omega_1$
on $M$ is called relatively exact modulo $\F$ if the restriction
of $\omega_1$ to each leaf $\L$ of $\F$ is exact i.e., there is a
meromorphic
 function $f$ on $\L$ such that $\omega_1\mid_{L}=df$.
\end{defi}
In fact we are interested in the above definition when $\F$ is
integrable. It is easy to check that a meromorphic 1-form
$\omega_1$ is relatively exact modulo $\F$ if and only if
\begin{equation}
\label{int}
\int_\delta \omega_1 =0
\end{equation}
for all closed cycles in the leaves of $\F$, where this integral is well-defined.

Relatively exact 1-forms have been studied by many authors.
Ilyashenko in \cite{ily} proves that if the integral of a
polynomial 1-form of degree $n$ along a continuous family of
level lines of a Morse type polynomial of two variables and of
degree $n$ is identically equal to zero, then the form is
relatively exact modulo the polynomial and then he proves that it
must be exact. This is generalized to higher dimensions in
\cite{pus}. The classification of relatively exact polynomial
1-forms modulo an arbitrary polynomial is done by P. Bonnet in
\cite{bon} and L. Gavrilov in \cite{gav} . J. Muci\~no in
\cite{muc} has classified a certain class
 of relatively exact 1-forms modulo  a Lefschetz pencil.
In a generalization of Ilyashenko's result to integrable foliations in $M$,
I  had to  classify another types of relatively exact 1-forms in \cite{hos1}.

Let $S$ be a complex curve. The meromorphic function
$f:M\rightarrow S$ is called non-composite if a general compact
$f$-fiber is irreducible. It is easy to see that $f:M\rightarrow
S$ is non-composite if and only if $f$ cannot
 be factored as a composite
\begin{equation}
\label{2may}
M\stackrel{f'}{\rightarrow }S'\stackrel {i}{\rightarrow} S
\end{equation}
where $S'$ is a complex curve and $i$ is a holomorphic map of degree greater than one.
\\
Let $\F$ be an integrable foliation. There exists a complex curve $S$ and a non-composite meromorphic function $f:M\rightarrow S$ such that $f$  is constant in the leaves of $\F$ (see \cite{gom}). We say that $f$ is a non-composite first integral of $\F$.
\\
Let $f:M\rightarrow S$ be non-composite as above. Any other first integral $f':M\rightarrow S'$ of $\F$ is factored as
$$
M\stackrel{f}{\rightarrow }S\stackrel {i}{\rightarrow} S'
$$
 where $i$ is a holomorphic function.
\\
Poincar\'e in his article \cite{poi} has studied integrable foliations in $\proj$ and has proved (p. 52-53) that every integrable foliation $\F$ in $\proj$ has a non-composite
first integral $M\rightarrow \pl$ (this is also a consequence of Stein factorization theorem).

Let $\F$ be an integrable foliation with the non-composite first integral $f:M\rightarrow S$ and $\omega_1$ a meromorphic 1-form in $M$ with the pole divisor
$$
D=\sum_{i=1}^{n} n_iD_i
$$
In what follows when we say that a meromorphic object (function, 1-form, ...) $Z$ has the pole divisor $D$, we mean that
$$
pol(Z)+D\geq 0
$$
Let
$$
D_0=\sum_{i=1}^{n} D_i
$$
be the reduced part of $D$.
\begin{theo}
Keeping the notations used above, assume that
\label{principal}
\begin{enumerate}
\item
All $f$-fibers are connected (for us a $f$-fiber does not contain
the indeterminacy points of $f$);
\item
$D$ is $\F$-invariant;
\item
Every component of $D$  is an irreducible $f$-fiber ($f$ may have
multiplicity along some $D_i$);
\item
There exists a  non $\F$-invariant Riemann surface embedded in $M-\RR$.
\end{enumerate}
Then every  relatively exact
meromorphic 1-form $\omega_1$ modulo $\F$ with the pole divisor $D$ has the form
\begin{equation}
\label{form1}
\omega_1=dg+\omega
\end{equation}
where $g$ is a meromorphic function on $M$ with the pole divisor
$D$ and $\omega$ is a meromorphic 1-form in $M$ inducing the
foliation $\F$ and with the pole divisor $D+D_0$.
\end{theo}

The above theorem is no more true if we assume that some
$f$-fibers are disconnected. P. Bonnet in \cite{bon} gives the
example $f=x(1+xy)$ in $\C^2$ having the disconnected fiber
$f=0$. The   1-forms $y^{k+1}dx+xy^kdy,\ k>0$ are relatively
exact modulo $f$ but they are not of the form ~(\ref{form1}).
Without the hypothesis of connectedness of $f$-fibers the
classification of relatively exact polynomial 1-forms modulo a
polynomial is done in \cite{bon}. In \cite{hos1} we have
classified relatively exact 1-form in $\proj$  with non-invariant
divisors. The forth condition is trivial for an algebraic
manifold with a Lefschetz pencil in it. I do not know yet that
the third and forth conditions in the above theorem is really
necessary or not. The above classification of relatively exact
1-form will suffice to our purpose.
\\
\pr
 Let $\RR$ be the indeterminacy set of $f$ and $C$ a  non
$\F$-invariant Riemann surface in $U=M-\RR$. For any point $x\in
U$ let
$$
L_x \cap C=\{p_1,p_2,\ldots,p_e\}
$$
where $L_x=f^{-1}(g(x))$ is the fiber through $x$ and $e$ is the
intersection number of $C$ with  a generic $f$-fiber ($p_i$'s are
counted with multiplicity). Define
$$
g: M\backslash(\cup_{i=1}^{n} D_i)\rightarrow \C
$$
$$
g(x)=\frac{1}{e}(\sum \int_{x}^{p_i}\omega _1)
$$
where $\int_{x}^{p_i}$ is an integral over a path in $L_x$ which
connects $x$ to $p_i$. Since the $f$-fibers are connected,
$\int_{x}^{p_i}$ is well-defined. The idea of the definition of
$g$ comes from the paper \cite{muc}. According to the hypothesis
this integral does not depend on the choice of the path
connecting $x$ to $p_i$ on $L_x$. The function $g$ is a
well-defined holomorphic function in $M\backslash(\cup_{i=1}^{n}
D_i)$. We claim that $g$ is a meromorphic function on $M$ with
pole divisor $D$. According to Levi extension theorem it is
enough to prove that $g$ is meromorphic in $U$.
\\
For instance let us prove that $g$ has a pole of order at most $n_i$ at
$U\cap D_i$. Let $m_i$ be the multiplicity  of $f$ along $D_i$, $c_i=f(D_i)$
be the value associated to $D_i$ and $x\in U\cap D_i$.
In a small neighborhood of the path connecting $x$ to $p_i$ the
function $(f-c_i)^{\frac{n_i}{m_i}}$ is a univalued holomorphic function and
$$
\int_{x}^{p_i} \omega_1=
(f-c_i)^{-\frac{n_i}{m_i}}\int_{x}^{p_i}(f-c_i)^{\frac{n_i}{m_i}}\omega_1
 $$
$(f-c_i)^{\frac{n_i}{m_i}}\omega_1$ is a holomorphic 1-form along
$U\cap D_i$ and therefore integrals of this type has a pole of
order at most $n_i$ at $\{D_i=0\}$. This implies that $g$ has
pole of order at most $n_i$ at $\{D_i=0\}$. Note that in the
above formulas we have chosen a local chart $z(c),c\in(S,c_i)$
 around $c_i$ and instead of $z\circ f-z(c_i)$ we have simply written $f-c_i$.
\\
Every integral $\int_{x}^{p_i}\omega _1$ satisfies the equation
$$
d(\int_{x}^{p_i}\omega _1)\wedge \omega_0=\omega_1\wedge \omega_0
$$
where $\omega_0$ is a meromorphic 1-form inducing $\F$. The function
$g$ is the mean value of these integrals and so
$$
dg\wedge \omega_0=\omega_1\wedge \omega_0 \Rightarrow (\omega_1-dg)\wedge \omega_0=0
$$
$$
\Rightarrow \omega_1=dg+\omega
$$
where $\omega$ is a meromorphic 1-form inducing $\F$ and with the 
pole divisor  $D+D_0$.
\qed

Let $\F$ be a holomorphic foliation in $M$ considered in Theorem
~\ref{principal}. Assume that  $\F$ has a non composite
meromorphic first integral $f$ ($S=\pl$). Denote by $D$ a generic
fiber of $f$.
\begin{cor}
\label{gij}
Every relatively exact 1-form $\omega_1$ in $M$ with the pole divisor $nD$ has the form
$$
\omega_1=dg+pdf
$$
where $g$ and $p$ are meromorphic functions with the pole divisors $nD$ and $(n-1)D$, respectively.
\end{cor}
\section{Foliations with a center in $\proj$}
\label{hargez}
\def\In{{\cal I}(a,b)}
\def\PP{{\cal P}_{a+1} \times {\cal P}_{b+1}}
\def\RR{\{F=0\}\cap\{G=0\}}
\def\INT{pGdF-qFdG}
\def\BBB{\backslash}
\def\Inm{{\cal I} _m(a,b)}
\def\t{\tilde}
Let $\F$ be a germ of singular foliation at $(\C ^2,0)$. We say that $0\in \C^2$ is a
center singularity of $\F$ or simply a center of $\F$, if there exists a germ of
holomorphic function
$f:(\C ^2,0)\rightarrow (\C ,0)$ which has non-degenerate critical point at $0\in\C^2$,
and the leaves of $\F$ near $0$ are given by $f=const.$. The point $0$ is also called a
Morse singularity of $f$. Morse lemma in the complex case implies that there exists a local
 coordinate
 system $(x,y)$ in $(\C^2,0)$  with $x(0)=0,\ y(0)=0$ and such that
 $f(x,y)=x^2+y^2$. Near the center the leaves of $\F$ are homeomorphic to a
  cylinder, therefore each leaf has a nontrivial closed cycle that will be
  called the Lefschetz vanishing cycle.

Let $\Sff$ be the space of degree $d$ holomorphic foliations in $\proj$ and 
$\M$ the closure of the set of foliations of degree $d$ and
 with at least one center in
$\Sff$. The following
example gives us a huge number of these foliations: Let $\tau
:\proj \rightarrow M$ be a holomorphic map between $\proj$ and a
complex compact manifold $M$ with $dim(M)\geq 3$. Let also $\F$
be a holomorphic codimension one singular foliation in $M$ (see
\cite{alci}). We say that $\tau$ has a tangency point $a\in \proj$
of order two with $\F$ if $\tau(a)$ is a regular point of the
foliation $\F$, $\F$ in a
  coordinate $(x,y)\in (\C^{n-1},0)\times(\C,0)$ around $a$ is given by
  $y=const.$ and $a$ is  a non-degenerate critical point of $y\circ\tau$.
  This says that the pullback foliation $\tau^{*}(\F)$ has a center at the point $a$.

I have learned the statement and proof of the following
proposition from A. Lins Neto.
\begin{pro}
\label{alge} $\M$ is an algebraic subset of $\Sff$.
\end{pro}
The proof of this proposition is given in the end of this section. Let  $\P _d$ be the set of polynomials of maximum degree $d$ in $\C ^2$ and
\[
(F,G)\in  \PP,\ \frac{a+1}{b+1}=\frac{q}{p},\ g.c.d.(p,q)=1
\]
The foliation $\F=\F(\INT)$ has the first integral:
\[
f:\proj\BBB (\RR) \rightarrow \S ,\
f(x,y)=\frac{F(x,y)^p}{G(x,y)^q}
\]
i.e., the leaves of the foliation $\F$ are contained in the level
surfaces of $f$. Let $\In$ be the closure of the set of the mentioned
holomorphic foliations in $\F(d)$. 
 Our main result in this section is the following:
\begin{theo}
\label{irreducible} $\In, a+b>2$ is an irreducible component of
$\M$, where $d=a+b$.
\end{theo}
This result is announced in \cite{hos1}.
We can restate our main theorem as follows: Let $\F\in\In$, $p$
one of the center singularities of $\F$ and $\F_t$ a holomorphic
deformation of $\F$ in $\Sff$, where $d=a+b$, such that its unique
singularity $p_t$ near $p$ is still a center.
\begin{theo}
\label{main}
 In the above situation, if $a+b>2$ then there exists an open dense subset $U$ of $\In$, such
 that for all $\F(\INT)\in U$, $\F_t$ admits a meromorphic first integral.
 More precisely, there
 exist polynomials $F_t$ and $G_t$ such that $\F_t=\F(pG_tdF_t-qF_tdG_t)$,
 where $F_t$ and $G_t$ are holomorphic in $t$ and $F_0=F$ and $G_0=G$.
\end{theo}
This theorem also says that the persistence of one center implies
the persistence of all other centers and dicritical singularities 
(the points of $\{F=0\}\cap\{G=0\}$).

The classification of degree two polynomial differential equations
was done by Dulac in \cite{dul}. Going to the language of
holomorphic foliations in $\proj$, instead of using the language
of polynomial differential equations,  this classification was
completed in \cite{alci1} for degree two holomorphic foliation in
$\proj$. This classification for some degree three differential
equations is done in \cite{llo}. Deformation of real Hamiltonian
equations with a center singularity, generating limit cycles from
the Lefschetz vanishing cycles of the center, has been one of the
methods of approach to Hilbert sixteen problem, for this  see
Roussarie's book \cite{rou} and its reference. Yu.S. Ilyashenko
in \cite{ily} shows that the persistence of a center after
deformation of a generic Hamiltonian equation implies that the
deformed equation is also Hamiltonian. He uses this fact to get a
certain number of limit cycles after deformation of Hamiltonian
equations. Theorem ~\ref{irreducible} is a generalization of the
Ilyashenko's result. The case $p=q=1$ was studied by J. Muci\~ no
in \cite{muc}. With more strong hypotheses on the tangent vector
of the deformation, i.e. the vector $\omega_1$ in
$\F_t(df+t\omega_1+h.o.t.)$, he succeeded in proving that
$\omega_1$ is a relatively exact 1-form and then classifies them.
This result, however, does not imply that the persistence of just
one center implies the persistence of all others.

{\bf Proof of Theorem ~\ref{irreducible}:}
Since $\In$ is
parameterized by $\PP$, $\In$ is an irreducible variety.
Let $\Inm$ be the subset of $\In$ containing $\F(\INT)$ with the
properties 1,2 before Theorem ~\ref{argentina}. $\Inm$ is an open dense
subset of $\In$ (see \cite{hos2}. Let $\F(\INT)\in\Inm$ and $\F_\epsilon(\INT+\epsilon\omega_1+h.o.t.)\in\Sff$ a deformation of $\F$ such that there is $\F_{\epsilon_i},\ i=1,2.\ldots$ with this property that $\F_{\epsilon_i}$ has a center $p_i$ near 
a fixed center $p$ of $\F$. Since $\F_\epsilon$ is of degree $d$,
$\frac{\omega}{FG}$ is a meromorphic 1-form in $\proj$ with poles only in $\{F=0\}$
and $\{G=0\}$ (there is no pole along the line at infinity).
The first Melnikov function of this deformation along
a Lefschetz vanishing cycle around $p$ is zero. Therefore
\begin{equation}
\label{inte} \int _{\delta} \frac{\omega _1}{FG} =0
\end{equation}
for all vanishing cycles $\delta$ in the leaves around  $p$ and by
theorem ~\ref{argentina}  we conclude that the equality
~(\ref{inte}) holds for all closed cycles $\delta$  in the leaves
of  the foliation $\F$, where the integral is defined. Partially
we obtain that the residue of $\frac{\omega _1}{FG}$ on a leaf
around any $q\in\{F=0\}\cap\{G=0\}$ is zero. Therefore the 1-form
$\frac{\omega _1}{FG}$ is relatively exact modulo
the foliation $\F(\INT)$. By Theorem ~\ref{principal}, there exist
 polynomials $B$ and $R$ of degree at most $d+2$ such that 
 \[
 \frac{\omega _1}{FG}=d(\frac{B}{FG})- \frac{R}{FG}(\frac{\INT}{FG})\Rightarrow
\]
\begin{equation}
\label{tan}
 \omega _1= \frac{FGdB-Bd(FG)-R(\INT )}{FG}
\end{equation}
This implies that
\[
F\mid B+pR,\ G\mid B-qR \Rightarrow
\]
\[
B+pR=-(p+q)FQ,\ B-qR=(p+q)GP\ \Rightarrow
\]
\[
B=pGP-qFQ,\ R=-GP-FQ
\]
where $P$ and $Q$ are two polynomials of degree at most $deg(F)$
and $deg(G)$, respectively. Substituting this in ~(\ref{tan}) we have
 \begin{equation}
\omega _1=pGdP-qPdG +pQdF-qFdQ 
\label{tang}
\end{equation}
Until now we have proved that 
\[
    [\omega _1] \in T_{\F}\M \Rightarrow
\]
\[
 \omega _1=pGdP-qPdG+pQdF-qFdQ,\ (P,Q)\in\PP\Rightarrow
\]
\[
[\omega_1]\in T_{\F}\In
\]
this and the fact that $\In\subset\M$ imply that
\[
T_{\F}\M =T_{\F}\In,\ \forall \F\in\Inm
\]
Since $\Inm$ is an open dense subset of $\In$, we conclude that
 $\In$ is an irreducible
component of $\M$. \qed

In the language of differential equations (see \cite{rou}) what
we have proved in Theorem ~\ref{irreducible} is the following:
The Bautin variety (zeros of Bautin ideal) of a Lefschetz
vanishing cycle of a generic $\F\in\In$ with $\Sff$, $d=a+b$, as
the deformation space, is $\In$. Note that the Bautin variety is
defined locally but $\In$ is a closed algebraic set in $\Sf$.

Let $X$ be an irreducible component of $\M$, $\F(\omega_0)\in X$ and $p$ be a
center of $\F$. There is a coordinate $(x,y)$ in a small neighborhood $U$ of $p$ such
 that in this coordinate $p=(0,0)$ and
\[
\omega_0=gd(\tilde{f}),\ \tilde{f}=\frac{1}{2}(x^2+y^2)+h.o.t.,\ g(0)\not =0
\]
where $\tilde{f}$ and $g$ are holomorphic functions on $U$.
Define $T^*_{\F}X$ as the set of all 1-forms
$[\omega_1]\in T_{\F}\Sf$ such that
\[
 \int_{\delta}\frac{\omega_1}{g}=0
\]
for all Lefschetz vanishing cycles in the leaves of $\F$ around
$p$.
By using the formula of first Melnikov function we know that
\begin{equation}
\label{strange}
T_{\F}X\subset T^*_{\F}X
\end{equation}
$X$ is called a good irreducible component of $\M$ if for a generic choice of
$\F\in X$, the equality holds in ~(\ref{strange}).
\\
The above definition does not depend on the choice of $g$. More
precisely, if
\[
\omega_0=gd\t{f}=g'd\t{f}'\Rightarrow \frac{g}{g'}d\t{f}=d\t{f}'\Rightarrow
d(\frac{g}{g'})\wedge \omega_0=0
\]
This means that $\frac{g}{g'}$ is constant on the leaves of $(\F, p)$ and so
\[
\int_{\delta}\frac{\omega_1}{g'}=\int_{\delta}\frac{g}{g'}\frac{\omega_1}{g}=
\frac{g}{g'}\int_{\delta}\frac{\omega_1}{g'}=0
\]
which implies that the definition does not depend on $g$.
\begin{pro}
$\In$ is a good irreducible component of $\M$.
\end{pro}
\pr The proof is the same as the proof of Theorem \ref{irreducible}.\qed
\\
{\bf Challenge:}
We have proved that for $d>2$
\[
{\cal I}(0,d),{\cal I}(1,d-1),{\cal I}(2,d-2), \ldots, {\cal
I}([\frac{d}{2}], d-[\frac{d}{2}])
\]
are good irreducible components of $\M$ (Note that $\In={\cal
I}(b,a)$). Which are the other irreducible components of $\M$? Is
any other irreducible component of $M$ good? For the cases
$d=1,2$ the complete description of the irreducible components of
$\M$ is done in \cite{dul} and \cite{alci1}. Let us introduce an
algebraic set in $\Sf$ which seems to be an irreducible component
of $\M$:
\\
Let $\F$ be a foliation in $\proj$ given by the polynomial 1-form
\begin{equation}
\label{log} \omega( f,\lambda)=\omega(f_1,\ldots,
f_r,\lambda_1,\ldots,\lambda_r)= f_1\cdots f_r \sum_{i=1}^r
\lambda_i \frac{df_i}{f_i}
\end{equation}
where the $f_i$'s are irreducible polynomials in $\C^2$, $\sum
d_i\lambda_i=0$ and $d_i=deg(f_i)$.
 $\F$ is called a logarithmic foliation and it has the
multi-valued first integral $f=f_1^{\lambda_1}\cdots
f_r^{\lambda_r}$ in $U=\proj\backslash (\cup_{i=1}^r \{f_i=0\})$.
We can prove that generically, the degree of $\F$ is
$d=\sum_{i=1}^r d_i-2$ and has  $d^2+d+1-\sum_{1\leq i<j\leq r}
d_id_j$ centers.

Let ${\cal L}(d_1,d_2,\ldots,d_r)$ be the set of all logarithmic
foliations of the above type. Is  ${\cal L}(d_1,d_2,\ldots,d_r)$
an ( a good) irreducible component of $\M$, where $d=\sum_{i=1}^r
d_i-2$? As far as I know, there is no any theorem similar to
~\ref{argentina} for logarithmic foliations.

Recently some developments in the above problem is made by the author
of these lines in \cite{hos4} for logarithmic differential equations in 
$\C^2$ (holomorphic foliations in $\proj$ with the line at infinity invariant). 

{\bf Proof of Proposition ~\ref{alge}:}
Let $\Mz$ be the set of all foliations in $\M$ with a center
at the origin $(0,0)\in\C^2\subset\proj$ and with a local first
integral of the type
\begin{equation}
\label{fff} f=xy+f_3+f_4+\cdots+f_n+h.o.t.
\end{equation}
in a neighborhood of $(0,0)$. Let us prove that $\Mz$ is an
algebraic subset of $\Sf$.
\\
Let $\F (\omega)\in\Mz$ and $\omega =\omega _1 +\omega _2 +\omega
_3 +\ldots +\omega _{d+1}$ be the homogeneous decomposition of
$\omega$, then in a neighborhood around $(0,0)$ in $\C^2$, we have
\[
\omega \wedge df =0 \Rightarrow (\omega _1 +\omega _2 +\omega _3
+\cdots +\omega _{d+1}) \wedge(d(xy) +df_3+df_4+\cdots)=0
\]
Putting the homogeneous parts of the above equation equal to
zero, we obtain
\begin{equation}
    \left\{ \begin{array}{ll}
        \omega _1 \wedge d(xy)=0 \Rightarrow \omega _1 =k.d(xy),\hbox{ $k$ is constant}  \\
        \omega _1 \wedge df_3 =-\omega _2 \wedge d(xy) \\
        \ldots   \\
        \omega _1 \wedge df_n =-\omega _2 \wedge df_{n-1} -\cdots-\omega
        _{n-1} \wedge d(xy) \\
        \ldots
    \label{ome}
       \end{array}
    \right.
\end{equation}
Dividing the 1-form $\omega$ by $k$, we can assume that $k=1$.
Let $\P_n$ denote the set of homogeneous polynomials of degree
$n$. Define the operator :
\[
    S_n:\P _n \rightarrow (\P _ndx \wedge dy)
\]
\[
 S_n(g)=\omega_1 \wedge d(g)
\]
We have
\[
S_{i+j}(x^iy^j)=d(xy)\wedge
d(x^iy^j)=(xdy+ydx)\wedge(x^{i-1}y^{j-1}(jxdy+iydx))
\]
\[
=(j-i)x^iy^j dx \wedge dy
\]
This implies that when $n$ is odd $S_n$ is bijective and so in
~(\ref{ome}), $f_n$ is uniquely defined  by the terms
$f_m,\omega_m$'s $m < n$,  and when $n$ is even
\[
Im(S_n)=A_ndx\wedge dy
\]
where $A_n$ is the subspace generated by the monomials $x^iy^j,\
i\not = j$. When n is even the existence of $f_n$ implies that
the coefficient of $(xy)^{\frac{n}{2}}$ in
\[
-\omega _2 \wedge df_{n-1} -\cdots-\omega_{n-1}\wedge d(xy)
\]
which is a polynomial,  say $P_n$, with variables
\[
\hbox{ coefficients of } \omega _2 \ldots \omega _{n-1},f_2
,\ldots ,f_{n-1}
\]
is zero.The coefficients of $f_i, i\leq n-1$ is recursively given
as polynomials in coefficients of $\omega_i,i\leq n-1$ and so the
algebraic set
\[
 X : P_4=0 \ \& \ P_6=0\ \&\ldots\&\ P_n=0\ \ldots
\]
consists of all foliations $\F$ in $\Sf$ which have a formal
first integral of the type ~\ref{fff} at $(0,0)$. From theorem A
in \cite{mat}, it follows that $\F$ has a holomorphic first
integral of the type ~(\ref{fff}). This implies that $\Mz=X$ is
algebraic. Note that by Hilbert zeroes theorem, a finite number
of $P_i$'s defines $\Mz$. The set $\M$ is obtained by the action
of the group of automorphisms of $\proj$ on $\Mz$. Since this
group is compact we conclude that $\M$ is  also algebraic. \qed

\section{Higher Melnikov functions} 
\label{melnikov}
In this section we follow the notations introduced in the first
lines of Section ~\ref{akharinseminar}. Let us suppose that $f$ is
non composite and  $D=M_\infty$ is a generic compact $f$-fiber.
The integrating factor $s$ of $\omega_0$ with
$$
\frac{\omega_0}{s}=df
$$
has $2D$ as the zero divisor.
 For simplicity let us write $\omega_\epsilon$ instead of $\frac{\omega_\epsilon}{s}$ and
 $\omega_i$ instead of $\frac{\omega_i}{s},\
 i=1,2,\ldots$.

\begin{theo}
Suppose that $\delta$ is a simple cycle defined in Section ~\ref{1feb01}. If
$M_1\equiv M_2\equiv\cdots\equiv M_{k}\equiv 0$ then
$$
M_{k+1}(t)=-\int_{\delta_t}(\sum_{i=1}^{k}p_i\omega_{k+1-i}+\omega_{k+1})
$$
where $p_i$ and $g_i$ are meromorphic functions in $M$ with the pole divisors $iD$ and $(i+1)D$, respectively,  and are
 defined recursively by
$$
\omega_i+p_idf+dg_i=-\sum_{j=1}^{i-1}p_j\omega_{i-j},\ i=1,2\ldots,k
$$
 \end{theo}
\pr The proof essentially follows \cite{rou} Proposition 6 p. 73.
We prove by induction on $k$. The case $k=1$ is proved in
Proposition ~\ref{1febb}. Let us suppose that the theorem is true
for $k-1$, i.e. if $M_1\equiv M_2\equiv\cdots\equiv M_{k-1}\equiv
0$ then
$$
M_{k}(t)=-\int_{\delta_t}(\sum_{i=1}^{k-1}p_i\omega_{k-i}+\omega_{k})
$$
Now suppose that $M_k\equiv 0$. Since $\delta$ is a simple cycle, the 1-form
\begin{equation}
\label{ach}
-(\sum_{i=1}^{k-1}p_i\omega_{k-i}+\omega_{k})
\end{equation}
is a relatively exact 1-form with the pole divisor $(k+1)D$ and so
by corollary ~\ref{gij} there exist $p_k$  and $g_k$ with the pole divisors $kD$ and $(k+1)D$, respectively, such that
$$
-(\sum_{i=1}^{k-1}p_i\omega_{k-i}+\omega_{k})=dg_k+p_kdf
$$
A direct expansion gives
$$
(1+\sum_{i=1}^{k}p_i\epsilon^i)\omega_\epsilon=d(f-\sum_{i=1}^{k}g_i\epsilon^i)+
(\sum_{i=1}^{k}p_i\omega_{k+1-i}+\omega_{k+1})\epsilon^{k+1}+O(\epsilon^{k+2})
$$
Let  $\delta_{t,h_\epsilon(t)}$ be a path in the leaf of
$\F_\epsilon$ through $t$ which connects $t$ to $h_\epsilon(t)$
along the path $\delta$. Since $\Sigma$ is parameterized by
$t=f\mid_{\Sigma}$, integrating the above equality over the path
$\delta_{t,h_\epsilon(t)}$ we have
$$
(h_\epsilon(t)-t)-(\sum_{i=1}^{k}g_i\epsilon^i)
\mid_{t}^{h_\epsilon(t)}+\epsilon^{k+1}
\int_{\delta_{t,h_\epsilon(t)}}(\sum_{i=1}^{k}p_i\omega_{k+1-i}+
\omega_{k+1})+O(\epsilon^{k+2})=0
$$
$\int_{\delta_{t,h_\epsilon(t)}}=\int_{\delta_t}+O(\epsilon)$ and
so by putting zero the coefficient of $\epsilon^{k+1}$ in the
above formula we get the desired equality.\qed

\smallskip
\leftline{Hossein Movasati}
\leftline{Institute for studies in
theoretical Physics and Mathematics, IPM} \leftline{School of
Mathematics} \leftline{P.O.Box : 19395-5746}
\leftline{Tehran-Iran} 
\leftline{Instituto de Matem\'atica
Pura e Aplicada, IMPA} \leftline{Estrada Dona Castorina, 110,
22460-320} \leftline{Rio de Janeiro, RJ, Brazil} \leftline{E-mails:
hossein@impa.br, movasati@ipm.ir}
\end{document}